\documentclass[preprint,3p,sort]{elsarticle}

\usepackage{amsfonts}
\usepackage{amsthm}
\usepackage{graphics}
\usepackage{epstopdf}
\usepackage{algorithmic}
\usepackage{ams math}
\usepackage{amssymb}
\usepackage{hyperref}
\usepackage[T1]{fontenc}
\usepackage{soul,color}
\usepackage{cases}
 \usepackage{lineno}
\newtheorem{theorem}{Theorem}[section]
\newtheorem{proposition}[theorem]{Proposition}
\newtheorem{lemma}[theorem]{Lemma}

\newtheorem{definition}{Definition}[section]
\newtheorem{remark}[theorem]{Remark}
\theoremstyle{example}
\newtheorem{example}[theorem]{Example}
\numberwithin{equation}{section}
\begin{document}


\begin{frontmatter}

%

\title{Nonlinear SBP operators for stable approximations of skew-symmetric hyperbolic initial boundary value problems}

\author[sweden,southafrica]{Jan Nordstr\"{o}m}
\cortext[secondcorrespondingauthor]{Corresponding author}
\ead{jan.nordstrom@liu.se}
\address[sweden]{Department of Mathematics, Applied Mathematics, Link\"{o}ping University, SE-581 83 Link\"{o}ping, Sweden}
\address[southafrica]{Department of Mathematics and Applied Mathematics, University of Johannesburg, P.O. Box 524, Auckland Park 2006, Johannesburg, South Africa}

\begin{abstract}
Summation-by-parts (SBP) operators in combination with simultaneous-approximation-terms (SATs) for boundary conditions provides a recipe that can be applied to well posed linear initial boundary value problems.
The SBP-SAT frame work guarantees semi-discrete energy stability and can be applied to most 
numerical discretization methods. 
High order accuracy 
comes without additional difficulty.

The extension of SBP-SAT to nonlinear initial boundary value problems has so far been somewhat unclear.  However, recent developments have shown how to construct stable nonlinear SAT terms. In this note we develop nonlinear SBP operators. These combined developments 
provide an SBP-SAT based recipe for stable approximations of nonlinear skew-symmetric hyperbolic initial boundary value problems.

\end{abstract}



\end{frontmatter}


\section{Introduction}
Energy and entropy stability for discretized well posed initial boundary value problems (IBVPs) require boundary procedures (including both the boundary conditions and the implementation techniques) at open boundaries for the continuous IBVPs that bound the surface terms in external data leading to a bounded energy rate. 
For linear hyperbolic IBVPs, the most common boundary  conditions are based on eigenvalue and eigenvector decomposition of the boundary term where one specifies the ingoing characteristic variables in terms  the outgoing ones and boundary data \cite{kreiss1970,Majda1975607,Higdon}.  Summation-by-parts (SBP) operators in combination with the simultaneous-approximation-term (SAT) technique can directly mimic the continuous analysis and lead to a bounded semi-discrete energy rate. The SBP-SAT frame work (see the reviews  \cite{svard2014review,fernandez2014review}) guarantees linear semi-discrete energy stability. SBP-SAT can be applied to most
numerical discretization methods, see e.g. the finite element \cite{abgrall2020analysis},  finite difference \cite{Carpenter1999341}, finite volume \cite{nordstrom2003finite} and spectral \cite{carpenter1996spectral} approximations.

A new type of provably stable characteristic boundary conditions for nonlinear  IBVPs was recently derived in  \cite{nordstrom2022linear-nonlinear,NORDSTROM2024_BC,Nordström2025_open,Nordström2026_LIN,Winters2026}. The main difference between the linear characteristic boundary conditions and the new nonlinear one is a scaling with the solution dependent eigenvalues for the nonlinear case. The scaling with eigenvalues originates from the nonlinear matrices in the IBVP and is included in the new nonlinear SAT formulation. In this note we will proceed in a similar manner and relate first derivative linear and nonlinear SBP operators to each other.  Analogous to the characteristic boundary conditions, the main difference between the linear and nonlinear case is that  the effect of the nonlinear matrices in the IBVP is accounted for.  The derivation is first done on the continuous side which clarify its dependence on skew-symmetric spatial operators.  In the second step we mimic the continuous derivation discretely to arrive at nonlinear first derivative SBP operators. 


\section{Continuous analysis}\label{sec:cont_theory}
Consider the following general nonlinear hyperbolic IBVP with both conservative and primitive terms
\begin{equation}\label{eq:nonlin_B}
\mathcal{M}  U_t + \partial_i(\mathcal{A}_i (U) U)+\mathcal{B}_i(U)\partial_iU=0,  \quad t \geq 0,  \quad  \vec x=(x_1,x_2,..,x_k) \in \Omega
\end{equation}
augmented with boundary conditions $B(U)-G(\vec x,t)=0$ at the boundary $\delta \Omega$ and an initial condition $U(\vec x,0)-F(\vec x)=0$ in $\Omega$. 
In (\ref{eq:nonlin_B}),  the Einstein summation convention is used, the subscript $t$ indicate differentiation in time and the operator $\partial_i$ differentiates in the $x_i$ spatial direction. Furthermore, $\mathcal{M}$ is a symmetric positive definite (or semi-definite) time-independent matrix that defines an energy norm (or semi-norm) $\|U\|^2_\mathcal{M} = \int_{\Omega} U^T \mathcal{M} U d\Omega$. 
The $n \times n$ matrices $\mathcal{A}_i,\mathcal{B}_i$ are smooth functions (each matrix element is smooth) of the $n$ component long smooth vector $U$, but otherwise arbitrary. 
\subsection{Energy conservation and boundedness}\label{energy_boundedness}
We start by investigating when the formulation (\ref{eq:nonlin_B}) leads to useful bounds.
\begin{definition}
The IBVP (\ref{eq:nonlin_B}) is energy conserving if $\|U\|^2_\mathcal{M}$ only changes due to boundary effects. It is energy bounded if $\|U\|^2_\mathcal{M}  < \infty$ as $t  \rightarrow  \infty$.
\end{definition}
\begin{proposition}\label{lemma:Matrixrelation}
The IBVP  (\ref{eq:nonlin_B}) is energy conserving if
\begin{equation}\label{eq:boundcond_B}
\mathcal{B}_i= \mathcal{A}_i^T, \quad i=1,2,..,k.
\end{equation}
It is energy bounded if it is energy conserving and the boundary conditions $B(U)-G=0$ are such that 
\begin{equation}\label{1Dprimalstab_B}
\oint\limits_{\partial\Omega}U^T  (n_i \mathcal{A}_i)   \\\ U \\\ ds = \oint\limits_{\partial\Omega} \frac{1}{2} U^T ((n_i \mathcal{A}_i)  +(n_i \mathcal{A}_i )^T) U \\\ ds \geq - \oint\limits_{\partial\Omega} G^TG \\\ ds .
\end{equation}
\end{proposition}
\begin{proof}
The energy method applied to (\ref{eq:nonlin_B}) yields
\begin{equation}\label{eq:boundaryPart1_B}
\frac{1}{2} \frac{d}{dt}\|U\|^2_M + \oint\limits_{\partial\Omega}U^T  (n_i \mathcal{A}_i)  \\\ U \\\ ds= \int\limits_{\Omega}(\partial_iU)^T  \mathcal{A}_i U - U^T \mathcal{B}_i (\partial_iU) \\\ d \Omega,
\end{equation}
where $(n_1,..,n_k)^T$ is the outward pointing unit normal. The right-hand side of (\ref{eq:boundaryPart1_B}) is cancelled by (\ref{eq:boundcond_B}) leading to energy conservation. If in addition (\ref{1Dprimalstab_B}) holds, time-integration leads to an energy bound.
\end{proof}
\begin{remark}
The relation (\ref{eq:boundcond_B})  leads to the nonlinear skew-symmetric formulation 
\begin{equation}\label{eq:nonlin_A}
\mathcal{M} U_t + (\partial_i\mathcal{A}_i+\mathcal{A}_i^T\partial_i)U=0.
\end{equation}
Proposition \ref{lemma:Matrixrelation}  shows that whatever form the original IBVP has, energy boundedness and conservation follows if it is rewritten as in (\ref{eq:nonlin_A}), see \cite{nordstrom2022linear-nonlinear, Nordstrom2022_Skew_Euler}.  Symmetric matrices as in \cite{Abarbanel19811, nordstrom_roadmap,oliger1978,nordstrom2005} are not required.
\end{remark}
\subsection{The nonlinear spatial operator}\label{energy_boundedness}
Consider the nonlinear skew-symmetric spatial operator 
\begin{equation}\label{eq:Spatial_Op_cont}
\mathcal{D}=\partial_i \mathcal{A}_i + \mathcal{A}_i^T \partial_i
\end{equation}
in (\ref{eq:nonlin_A}). 
We can prove the following Lemma.
\begin{lemma}\label{lemma:Integration_by_parts}
The nonlinear differential operator in  (\ref{eq:Spatial_Op_cont}) is an integration-by-parts (IBP) operator satisfying
\begin{equation}\label{eq:IBP_Part0}
\int\limits_{\Omega} V^T (\mathcal{D}U) + (\mathcal{D}V)^T U d \Omega =  \oint\limits_{\partial\Omega} n_i \left (V^T (\mathcal{A}_i U)  +(\mathcal{A}_i V)^T U  \right) \\\ ds= \oint\limits_{\partial\Omega}V^T  \left(n_i (\mathcal{A}_i+\mathcal{A}_i^T)\right)  \\\ U \\\ ds.
\end{equation}
\end{lemma}
\begin{proof}
IBP of the two elements in the utmost left term in (\ref{eq:IBP_Part0}) leads to 
\begin{equation}\label{eq:IBP_Part1}
\int\limits_{\Omega} V^T \mathcal{D} U d \Omega =  \oint\limits_{\partial\Omega}V^T  (n_i (\mathcal{A}_i+\mathcal{A}_i^T))  \\\ U \\\ ds - \int\limits_{\Omega} \left[(\partial_i V)^T  \mathcal{A}_i + ( \partial_i(\mathcal{A}_iV))^T  \right] U \\\ d \Omega.
\end{equation}
Next, we define the differentiation operator $\partial_i^T$ to be such that it operates to the left. The usual transpose rule then formally gives  $(\partial_i V)^T=V^T \partial_i^T$  and $ ( \partial_i (\mathcal{A}_iV))^T =V^T \mathcal{A}_i^T\partial_i^T$. Applying this to the square bracket in (\ref{eq:IBP_Part1}) and repeatedly using the transpose rule gives
\begin{equation}\label{eq:IBP_Part2}
 \left[(\partial_i V)^T  \mathcal{A}_i + ( \partial_i(\mathcal{A}_iV))^T  \right] U =  V^T \left[  \partial_i^T \mathcal{A}_i + \mathcal{A}_i^T\partial_i^T  \right] U = V^T \left[  \mathcal{A}_i^T \partial_i + \partial_i \mathcal{A}_i  \right]^T U  = (\mathcal{D} V)^T U.
\end{equation}
By combining (\ref{eq:IBP_Part1}) and (\ref{eq:IBP_Part2}) we obtain the IBP formula (\ref{eq:IBP_Part0}) which proves the Lemma. 
\end{proof}

\section{Semi-discrete analysis}\label{sec:disc_theory}
We will mimic the derivation in the previous Section using SBP operators on matrix form $D_i$ which directly mimic IBP in the discrete setting.
We start by considering the semi-discretized version of (\ref{eq:nonlin_A})
\begin{equation}\label{eq:disc_nonlin_B}
M \vec  U_t + D_i (A_i  \vec U)+A_i^T D_i \vec U=0,\quad  \text{or equivalently} \quad M \vec  U_t + D_V(\vec U)\vec U=0.
\end{equation}
In (\ref{eq:disc_nonlin_B}), the SBP operator on matrix form $D_i$ approximates the derivative operator $\partial_i$, the vector $\vec U$ approximates the vector $U$ in all nodes and the matrices $M, A_i$  approximates the matrices $\mathcal{M},\mathcal{A}_i$ respectively in all nodes. The matrices $A_i$ have diagonal elements and the  SBP operator  $D_i$ satisfies 
\begin{equation}\label{eq:SBP_def}
P_V D_i =Q_i, \quad Q_i+Q_i^T = E^T P_B N_i E.
\end{equation}
In (\ref{eq:SBP_def}), $E$ is a projection operator that projects the solution $\vec U$ from the volume nodes in $\Omega$ to the boundary nodes on $\partial\Omega$  \cite{Lundquist201849}. The diagonal matrix $N_i$ contain the normal vectors $n_i$ in each boundary node.  Furthermore, the diagonal 
boundary quadrature  matrix  $P_B$ define the $L_2$ norm $\| E \vec  U \|^2_{P_B}= (E \vec U)^T P_B (E \vec  U)$.

Following the continuous analysis  
we have also defined the nonlinear difference operator
\begin{equation}\label{eq:Spatial_Op_disc}
D_V=D_i A_i +  A_i^T D_i \quad  \text{factorized as} \quad P_V D_V=Q_V
\end{equation}
where  $P_V$ is a diagonal volume quadrature matrix defining the  $L_2$ norm $\| \vec  U \|^2_{P_V}=\vec  U^T P_V \vec  U$. The matrix $Q_V$ will be shown below to be almost skew-symmetric.  
With the definitions (\ref{eq:disc_nonlin_B}),(\ref{eq:SBP_def}),(\ref{eq:Spatial_Op_disc}) and the related information, we can now prove the following Proposition.
\begin{proposition}\label{lemma:summation_by_parts} 
The nonlinear difference operator in $D_V$ (\ref{eq:Spatial_Op_disc}) is an SBP operator satisfying
\begin{equation}\label{eq:IBP_Part0_disc}
\vec V^T P_V (D_V \vec  U) + (D_V \vec V) ^T P_V \vec  U =  (E \vec V)^T P_B N_i (E A_i \vec  U) + (E A_i \vec V)^T P_B N_i (E \vec  U).
\end{equation}
\end{proposition}
\begin{proof}
Combine the utmost left term in (\ref{eq:IBP_Part0_disc}, the first part of (\ref{eq:SBP_def}) and the fact that $P_VA_i=A_iP_V$ to get
\begin{equation}\label{eq:IBP_Part1_disc}
\vec V^T P_V (D \vec  U) =\vec V^T P_V (D_i A_i +  A_i^T D_i) \vec  U =  \vec V^T (Q_i A_i +  A_i^T Q_i) \vec  U =  \vec V^T Q_V \vec  U.
\end{equation}
The second part of (\ref{eq:SBP_def}) (corresponding to IBP) separate out the volume and boundary terms as
\begin{equation}\label{eq:IBP_Part2_disc}
\vec V^T (Q_i A_i +  A_i^T Q_i) \vec  U = - \vec V^T (Q_i^T  A_i  +  A_i^T Q_i^T)  \vec  U + (E  \vec V)^T P_B N_i (E A_i \vec  U)  + (E A_i \vec V)^T  P_B N_i (E \vec  U). 
\end{equation}
The numerically integrated boundary terms by $P_B$ in (\ref{eq:IBP_Part2_disc}) now corresponds to the middle continuous boundary integral in (\ref{eq:IBP_Part0}). The volume term can be rewritten with the aid of (\ref{eq:IBP_Part1_disc}) as
\begin{equation}\label{eq:IBP_Part3_disc}
- \vec V^T (Q_i^T  A_i  +  A_i^T Q_i^T)  \vec  U  = -  \vec V^T (A_i^T Q_i  +  Q_i A_i)^T  \vec  U = - \vec V^T Q_V^T \vec  U =- \vec V^T D^T P_V \vec  U = - (D  \vec V)^T P_V  \vec U
\end{equation}
which combined with (\ref{eq:IBP_Part1_disc}) and (\ref{eq:IBP_Part2_disc})  proves the proposition. 
\end{proof}

\section{Summary}\label{sec:conclusion}
We have developed nonlinear first derivative SBP operators which include the effect of solution dependent matrices. 
The derivation is first done on the continuous side which clarify its dependence on skew-symmetric spatial derivative operators.
In the second step we mimic the continuous derivation discretely to arrive at nonlinear first derivative SBP operators. This development together with the previous development of nonlinear SAT boundary conditions provides a direct extension of the SBP-SAT framework to nonlinear skew-symmetric hyperbolic initial boundary value problems.
\section*{Acknowledgment}
J.N. was supported by University of Johannesburg Global
Excellence and Stature Initiative Funding.


\bibliographystyle{elsarticle-num}
\bibliography{References_Jan}

@article{nordstrom2005,
Author = {Jan Nordstr\"{o}m and Magnus Sv\"{a}rd},
Journal = {SIAM J. Numer. Anal.},
Pages = {1231--1255},
Title = {Well posed boundary conditions for the {N}avier--{S}tokes equations},
Volume = {43},
Year = {2005},
}

@ARTICLE{Winters2026,
	author = {Winters, Andrew R. and Kopriva, David A. and Nordström, Jan},
	title = {Numerical boundary flux functions that give provable bounds for nonlinear initial boundary value problems with open boundaries},
	year = {2026},
	journal = {Journal of Computational Physics},
	volume = {559},
	type = {Article},
}

@ARTICLE{Majda1975607,
	author = {Majda, Andrew and Osher, Stanley},
	title = {Initial‐boundary value problems for hyperbolic equations with uniformly characteristic boundary},
	year = {1975},
	journal = {Communications on Pure and Applied Mathematics},
	volume = {28},
	number = {5},
	pages = {607 – 675},
	type = {Article},
}

@article{nordstrom2022linear-nonlinear,
author={Nordström,J.},
year={2022},
title={Nonlinear and linearised primal and dual initial boundary value problems: When are they bounded? How are they connected?},
journal={Journal of Computational Physics},
volume={455 111001},
language={English},
}

@ARTICLE{Nordstrom2022_Skew_Euler,
	author = {Nordström, Jan},
	title = {A skew-symmetric energy and entropy stable formulation of the compressible {E}uler equations},
	year = {2022},
	journal = {Journal of Computational Physics},
	volume = {470 111573},
	type = {Article},
}

@article{NORDSTROM2024_BC,
title = {Nonlinear boundary conditions for initial boundary value problems with applications in computational fluid dynamics},
journal = {Journal of Computational Physics},
volume = {498 112685},
year = {2024},
author = {Jan Nordström},
}

@ARTICLE{Nordström2025_open,
	author = {Nordström, Jan},
	title = {Open boundary conditions for nonlinear initial boundary value problems},
	year = {2025},
	journal = {Journal of Computational Physics},
	volume = {530},
	type = {Article},
	publication_stage = {Final},
	source = {Scopus},
}

@ARTICLE{Carpenter1999341,
	author = {Carpenter, Mark H. and Nordström, Jan and Gottlieb, David},
	title = {A Stable and Conservative Interface Treatment of Arbitrary Spatial Accuracy},
	year = {1999},
	journal = {Journal of Computational Physics},
	volume = {148},
	number = {2},
	pages = {341 – 365},
	type = {Article},
}

@article{Higdon,
author = {Higdon, Robert L.},
title = {Initial-Boundary Value Problems for Linear Hyperbolic System},
journal = {SIAM Review},
volume = {28},
number = {2},
pages = {177-217},
year = {1986},
}

@article{svard2014review,
  title={Review of summation-by-parts schemes for initial--boundary-value problems},
  author={Sv{\"a}rd, Magnus and Nordstr{\"o}m, Jan},
  journal={Journal of Computational Physics},
  volume={268},
  pages={17--38},
  year={2014},
}

@article{fernandez2014review,
  title={Review of summation-by-parts operators with simultaneous approximation terms for the numerical solution of partial differential equations},
  author={Fern{\'a}ndez, David C Del Rey and Hicken, Jason E and Zingg, David W},
  journal={Computers \& Fluids},
  volume={95},
  pages={171--196},
  year={2014},
}

@ARTICLE{Abarbanel19811,
author={Abarbanel, S. and Gottlieb, D.},
title={Optimal time splitting for two- and three-dimensional {N}avier-{S}tokes equations with mixed derivatives},
journal={Journal of Computational Physics},
year={1981},
volume={41},
number={1},
pages={1-33},
}

@article{kreiss1970,
	title={Initial boundary value problems for hyperbolic systems},
  	author={Kreiss, Heinz-Otto},
  	journal={Commun. Pur. Appl. Math.},
 	volume={23},
  	number={3},
  	pages={277--298},
  	year={1970}}

@ARTICLE{Lundquist201849,
	author = {Lundquist, Tomas and Malan, Arnaud and Nordström, Jan},
	title = {A hybrid framework for coupling arbitrary summation-by-parts schemes on general meshes},
	year = {2018},
	journal = {Journal of Computational Physics},
	volume = {362},
	pages = {49 – 68},
        type = {Article}
}

@article{nordstrom2003finite,
  title={Finite volume methods, unstructured meshes and strict stability for hyperbolic problems},
  author={Nordstr{\"o}m, Jan and Forsberg, Karl and Adamsson, Carl and Eliasson, Peter},
  journal={Applied Numerical Mathematics},
  volume={45},
  number={4},
  pages={453--473},
  year={2003},
  publisher={Elsevier}
}

@article{carpenter1996spectral,
  title={Spectral methods on arbitrary grids},
  author={Carpenter, Mark H and Gottlieb, David},
  journal={Journal of Computational Physics},
  volume={129},
  number={1},
  pages={74--86},
  year={1996},
  publisher={Elsevier}
}

@ARTICLE{Nordström2026_LIN,
	author = {Nordström, Jan},
	title = {Linear and nonlinear boundary conditions: What’s the difference?},
	year = {2026},
	journal = {Journal of Computational Physics},
	volume = {550},
    number = {114649},
	type = {Article},
}

@article{abgrall2020analysis,
  title={Analysis of the {SBP-SAT} stabilization for finite element methods part I: linear problems},
  author={Abgrall, R{\'e}mi and Nordstr{\"o}m, Jan and {\"O}ffner, Philipp and Tokareva, Svetlana},
  journal={Journal of Scientific Computing},
  volume={85},
  number={2},
  pages={1--29},
  year={2020},
  publisher={Springer}
}

@article{nordstrom_roadmap,
	Author = {Jan Nordstr\"{o}m},
	Journal = {{Journal of  Scientific  Computing}},
	Number = {1},
	Pages = {365--385},
	Title = {A roadmap to well posed and stable problems in computational physics},
	Volume = {71},
	Year = {2017}}

@article{oliger1978,
	Author = {Oliger, Joseph and Sundstr{\"o}m, Arne},
	Journal = {SIAM J. Appl. Math.},
	Number = {3},
	Pages = {419--446},
	Title = {Theoretical and practical aspects of some initial boundary value problems in fluid dynamics},
	Volume = {35},
	Year = {1978}}

\end{document}